\documentclass{article}
\usepackage{amssymb,amsmath,amsthm}

\def\={~=~}
\newcommand{\lp}{\left(}
\newcommand{\lb}{\left\lbrack}
\newcommand{\rp}{\right)}
\newcommand{\rb}{\right\rbrack}

\newtheorem{lemma}{Lemma}

\begin{document}

\title{A Class of  Continued Radicals}
\author{Costas J. Efthimiou}
\date{}
\maketitle

\begin{abstract}
We compute the limits of a class of continued radicals extending the results of a previous note in which only periodic radicals of the class were considered.
\end{abstract}

%_____________________________________________________________________________
\section{Introduction.}

In \cite{Efthimiou} the author  discussed the values for a class of periodic continued radicals of the form
{\footnotesize
\begin{equation}
         a_0\sqrt{2+a_1\sqrt{2+a_2\sqrt{2+a_3\sqrt{2+\cdots}}}}     ~,
\label{eq:OurRadical}
\end{equation}
}
where for some positive integer $n$,
$$
     a_{n+k} \= a_k~,~~~k=0,1,2,\dots~,
$$     
and
$$
      a_k\in\{-1,+1\}~,~~~k=0,1,\dots,n-1~.
$$
It was also shown that the radicals given by equation \eqref{eq:OurRadical} have limits two times the
fixed points of the Chebycheff polynomials $T_{2^n}(x)$, thus unveiling an interesting relation between these topics.

 In \cite{ZH}, the authors defined the set $S_2$ of all continued radicals of the form \eqref{eq:OurRadical} (with $a_0=1$) and they
  investigated some of their properties by assuming that  the limit of the radicals exists.  In particular, they showed that all elements 
  of $S_2$ lie between 0 and 2,  any two radicals cannot be equal to each other, and $S_2$ is uncountable.

My previous  note hence partially bridged this gap but left unanswered the question `\textit{what are the limits if the radicals 
are not periodic?}'  I answer the question in this note. The result is easy to establish, but I realized it only as I was reading the proof of 
my previous note. Such is the working of the mind!

%_____________________________________________________________________________
\section{The Limits.}

 Towards  the desired result, I present the following lemma from \cite{Shklarsky}, also used in the periodic case, which is an extension of the well known trigonometric formulas of the angles $\pi/2^n$.
\begin{lemma}
For $a_i\in\{-1,1\}$, with $i=0,1,\dots,n-1$, we have that
{\footnotesize
$$
  2\, \sin  \lb \lp a_0+{a_0a_1\over2}+\dots+{a_0a_1\cdots a_{n-1}\over2^{n-1}} \rp {\pi\over4} \rb
  \=
  a_0\sqrt{2+a_1\sqrt{2+a_2\sqrt{2+\dots+a_{n-1}\sqrt{2}}}}~.
$$ 
}
\end{lemma}
The  lemma is  proved in \cite{Shklarsky} using induction.

According to this lemma, the partial sums of  the continued radical  \eqref{eq:OurRadical} are given by
$$
   x_n \=  2\sin  \lb \lp a_0+{a_0a_1\over2}+\dots+{a_0a_1\cdots a_{n-1}\over2^{n-1}} \rp {\pi\over4} \rb~.
$$
The series
\begin{equation*}
   a_0+{a_0a_1\over2}+\dots+{a_0a_1\cdots a_{n-1}\over2^{n-1}} +\cdots
\end{equation*}
is absolutely convergent and thus it converges to some number $a$. Therefore, the original continued radical
 converges to the real number
$$
 x \= 2\sin{a\pi\over4}~.
$$
We can find a concise formula for $x$. For this calculation it is more useful to use the products
$$
    P_m \= \prod_{k=0}^{m} a_k~,~~~\text{for } m=0,1,2,\dots~,
$$
which take the values $\pm1$.  We will refer to these as partial parities. (When the pattern is periodic of period $n$ only the first $n$ 
parities  $P_0, P_1, \dots, P_{n-1}$ are independent.)  Using the notation with the partial parities, set
\begin{eqnarray*}
    a  &=&  P_0 + {P_1\over2} + {P_2\over 2^2} + \cdots+  {P_{n-1}\over 2^{n-1}} +   {P_n\over 2^{n}} +  \cdots  ~.
\end{eqnarray*}

We now define 
$$
     Q_m \= {1+P_m\over2}~.
$$
Since $P_m\in\{-1,1\}$, it follows that $Q_m\in\{0,1\}$. Inversely, $P_m=2Q_m-1$. Thus
\begin{equation*}
   a   \=   \sum_{m=0}^\infty  {P_m\over 2^m}
                 \=   \sum_{m=0}^\infty  {Q_m\over 2^{m-1}} - \sum_{m=0}^\infty  {1\over 2^m}
                 \=   4\, \sum_{m=0}^\infty  {Q_m\over 2^{m+1}} - 2~.
\end{equation*}
Notice that the sum 
$$
       Q \= \sum_{m=0}^\infty  {Q_m\over 2^{m+1}}
$$
 in the previous equation is the number $Q$  whose binary expression is $0.Q_0Q_1\cdots Q_{n-3}Q_{n-2}\cdots$. 
Therefore $a=4Q-2$.
 In \cite{ZH}, the authors noticed that all continued radicals of the form \eqref{eq:OurRadical} (with $a_0=1$) are in one-to-one 
 correspondence with the set  of decimals between 0 and 1 as written in binary notation (and that's how they determined that the set 
 $S_2$ is uncountable). But, with the above calculation, this correspondence is made deeper. It gives the limit of the radical 
 \eqref{eq:OurRadical} as follows
$$
        x \=  -2\cos \left( Q\pi \right) ~.
$$
For example, if $a_k=1$ for all $k$, then also $Q_k=1$ for all $k$  and the number $Q=0.111111111\cdots$ written in the binary system 
is the  number $Q=1$ in the decimal system; hence $x=2$. We thus recover the well known result
{\footnotesize
\begin{equation*}
           2 \= \sqrt{2+\sqrt{2+\sqrt{2+\sqrt{2+\cdots}}}}     ~.
\end{equation*}
}

%_____________________________________________________________________________
\section{Conclusion.}

Having found the limit of \eqref{eq:OurRadical}, the next obvious question is to determine the limit of the radical
{\footnotesize
\begin{equation*}
         a_0\sqrt{y+a_1\sqrt{y+a_2\sqrt{y+a_3\sqrt{y+\cdots}}}}     ~,
\end{equation*}
}
for values of the variable $y$ that make the radical (and the limit) well defined. However, a direct application of the above method fails and
so far a convenient variation has been elusive.  Therefore, the limit of the last radical in the general case remains an open problem although it is
known in at least two cases \cite{ZH}.

%_____________________________________________________________________________

%___________________________________________________________________________
\bigskip

\noindent\textit{
Department of Physics, University of Central Florida, Orlando, FL 32816 \\
costas@physics.ucf.edu
}

\end{document}